\documentclass[reqno]{amsart}
\usepackage[centertags]{amsmath}
\usepackage{graphicx}
\usepackage{amsfonts}
\usepackage{amssymb}
\usepackage{amsthm}
\usepackage[top=1in, bottom=1in, left=1in, right=1in]{geometry}

\usepackage[pdfborder={0 0 0},
   pdftitle={MODULUS OF UNBOUNDED VALENCE SUBDIVISION RULES},
  pdfauthor={BRIAN RUSHTON}, pdftex, bookmarks=true, bookmarksnumbered = true]{hyperref}

\theoremstyle{definition}

\theoremstyle{definition}

\theoremstyle{definition}

\theoremstyle{remark}

\theoremstyle{definition}

\newcommand{\parlengths}{\setlength{\parindent}{0pt}}
\begin{document}
\date{\today}

\pdfbookmark[1]{MODULUS OF UNBOUNDED VALENCE SUBDIVISION RULES}{user-title-page}

\title{Modulus of unbounded valence subdivision rules}

\author{Brian Rushton}
\address{Department of Mathematics, Brigham Young University, Provo, UT 84602, USA}
\email{lindianr@gmail.com}

\begin{abstract}
Cannon, Floyd and Parry have studied the modulus of finite subdivision rules extensively. We investigate the properties of the modulus of subdivision rules with linear and exponential growth at every vertex, using barycentric subdivision and a subdivision rule for the Borromean rings as examples. We show that the subdivision rule arising from the Borromean rings is conformal, and conjecture that the subdivision rules for all alternating links are conformal. We show that the 1,2,3-tile criterion of Cannon, Floyd, and Parry is sufficient to prove conformality for linear growth, but not exponential growth. We show that the criterion gives a weaker form of conformality for subdivision rules of exponential growth at each vertex. We contrast this with the known, bounded-valence case, and illustrate our results with circle packings using Ken Stephenson's Circlepack.
\end{abstract}

\maketitle\parlengths
\section{Introduction}\label{Introduction}

A subdivision rule is a recursive way of combinatorially dividing a tiling on a surface into a smaller, more refined tiling. Barycentric subdivision is an example, as is hexagonal refinement, both of which we will explore later. Subdivision rules have been studied extensively by Cannon, Floyd, and Parry in an attempt to prove that all hyperbolic groups with a 2-sphere at infinity are hyperbolic 3-manifold groups; as a corollary, this would give a simple proof of the hyperbolization conjecture, now proved by Perelman.

While investigating these hyperbolic groups with a 2-sphere at infinity, they showed (in \cite{hyperbolic}) that all such groups have an associated finite subdivision rule on the sphere. Cannon showed in \cite{conformal} that, if a subdivision rule associated to a group satisfies two simple axioms, then the group acts on the sphere by M\"{o}bius transformations, and on hyperbolic 3-space by isometries, cocompactly and properly discontinuously.

To state the axioms, we need to define combinatorial modulus. Combinatorial modulus is a direct analog of the modulus of a topological annulus or quadrilateral in complex analysis. In that setting, every topological annulus or quadrilateral is conformally equivalent to a standard annulus or quadrilateral, and the modulus is defined by how thick or thin the resulting annulus or quadrilateral is.

A conformal map can be thought of as a change in metric; finding the classical modulus, then, is finding an optimal metric in some sense. Combinatorial modulus is defined to mirror this.

A tiling $T$ of a ring $R$ (i.e., a closed annulus) gives two invariants, $M_{sup} (R,T)$ and $m_{inf} (R,T)$, called \textbf{approximate moduli}. These are similar to the classical modulus of a ring. They are defined by the use of \textbf{weight functions}. A weight function $\rho$ assigns a non-negative number called a \textbf{weight} to each tile of $T$. Every path in $R$ can be given a length, defined to be the sum of the weights of all tiles in the path. We define the \textbf{height} $H(\rho)$ of $R$ under $\rho$ to be the infimum of the length of all possible paths connecting the inner boundary of $R$ to the outer boundary. The \textbf{circumference} $C(\rho)$ of $R$ under $\rho$ is the infimum of the length of all possible paths circling the ring (i.e. not nullhomotopic in R). The \textbf{area} $A(\rho)$ of $R$ under $\rho$ is defined to be the sum of the squares of all weights in $R$. Then we define $M_{sup} (R,T)=\mathop{\sup}\limits_{\rho} \frac{H(\rho)^2}{A(\rho)}$ and $m_{inf} (R,T)=\mathop{\inf}\limits_{\rho} \frac{A(\rho)}{C(\rho)^2}$. Note that they are invariant under scaling of the metric.

A sequence $T_1,T_2,...$ of tilings is \textbf{conformal ($K$)} if mesh approaches 0 and:
\begin{enumerate}
\item for each ring $R$, the approximate moduli $M_{sup}(R,T_i)$ and $m_{inf}(R,T_i)$, for all $i$ sufficiently large, lie in a single interval of the form $[r,Kr]$; and
\item given a point $x$ in the surface, a neighborhood $N$ of $x$, and an integer $I$, there is a ring $R$ in $N\setminus\{x\}$ separating x from the complement of $N$, such that for all large $i$ the approximate moduli of $R$ are all greater than $I$.
\end{enumerate}

Note that mesh approaching 0 is independent of topological metric.

These are the two axioms mentioned earlier. The first axiom is similar to equicontinuity of maps, and the second shows that points do not blow up. In Section \ref{Barycentric}, we give an example of a subdivision rule satisfying axiom 1 but not axiom 2.

One way to achieve an intuitive understanding of these axioms is by circle packings. Ha\"{i}ssinsky showed that, if a sequence of tilings is conformal, then some subsequence of its circle packings will converge to a limit map (if we normalize the circle packings by fixing the image of three points). Thus, throughout this paper, we will use circle packings as a way of studying modulus.

\section{Barycentric Subdivision}\label{Barycentric}

Cannon and Swenson have shown \cite{hyperbolic} that the subdivision rules arising from closed hyperbolic manifolds are conformal. Later, Cannon, Floyd and Parry \cite{subdivision} showed that barycentric subdivision does not satisfy Axiom 1. We show that by altering the definition of weightings, we can make barycentric subdivision satisfy Axiom 1 on a basis of annuli, but not Axiom 2.

We will use the Layer Theorem \cite{subdivision}. This theorem says that the modulus of an annulus is greater than the sum of the moduli of any disjoint collection of essential sub-annuli. In other words, if you stack several rings together, you get a ring at least as big as all of them put together.

\begin{figure}
\begin{center}
\scalebox{1.0}{\includegraphics{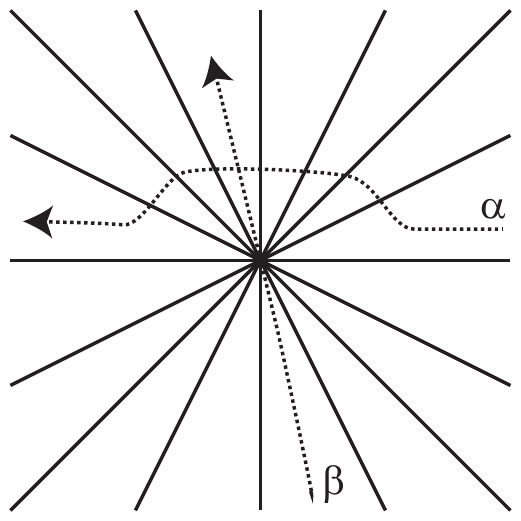}} \caption{Fat paths go around vertices, skinny paths go through them.} \label{FatAndSkinny}
\end{center}
\end{figure}

The modulus of a quadrilateral or annulus corresponds to an optimal weighting of the shinglings of the quadrilateral. The reason that unbounded valence subdivision rules have been difficult to study before is that the optimal weight vectors corresponding to a quadrilateral or annulus with tilings of large valence are asymmetrical. Specifically, Cannon, Floyd and Parry showed \cite{squaring} that optimizing a quadrilateral to measure its height in the traditional, `fat' way (counting all shingles that intersect a path) measures its width in a `skinny' way (counting any one set of shingles that covers the path). See Figure \ref{FatAndSkinny}. If many edges come into a single vertex, the fat height paths have to go all the way around the vertex, while the skinny width paths can sneak through the vertex. Thus, even rotationally symmetric tilings with high valence will have degenerate moduli.

\begin{figure}
\begin{center}
\scalebox{1.2}{\includegraphics{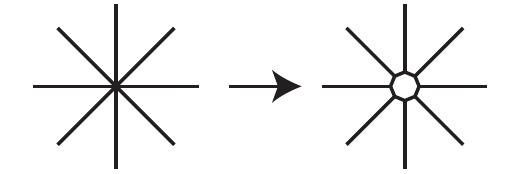}} \caption{Blowing up a vertex.} \label{BaryBlowup}
\end{center}
\end{figure}

My idea is to show that barycentric subdivision satisfies Axiom 1 on a basis of annuli if we change how we place weights. We either place weights at the vertices, or blow up the vertices by replacing each vertex with a closed disk (see Figure \ref{BaryBlowup}). Blowing up the vertices gives a new tiling with all vertices have valence 3. Placing weights at the vertices (and requiring paths to be edge paths) is equivalent to taking the dual tiling, where, again, all vertices have valence 3 (because all tiles of barycentric subdivision are triangles). In these settings, `skinny' and `fat' sets of shingles are the same, because any two tiles sharing a vertex share an edge. Thus, the height and width are measured in a symmetrical way. Because of this symmetry, Theorem 2.4.5.1 of \cite{squaring} shows that $m_{inf}=M_{sup}$ for all rings at all stages of subdivision. Thus, we can speak of `the' combinatorial modulus.

Our future calculations will not depend on the method we choose, as long as fat and skinny are the same. One motivation in putting weights on the vertices is circle packings; a circle packing assigns a weight to each vertex (i.e the radius of the circle).

\begin{figure}
\begin{center}
\scalebox{.8}{\includegraphics{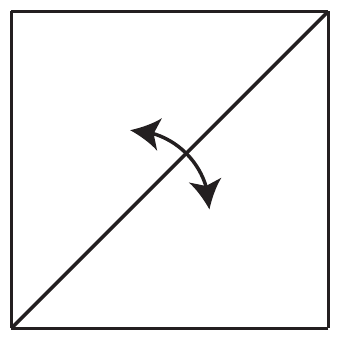}} \caption{A square in a triangulation.} \label{BaryReflect1}
\end{center}
\end{figure}

\begin{figure}
\begin{center}
\scalebox{.8}{\includegraphics{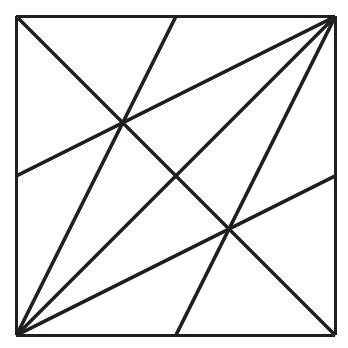}} \caption{The same square after one subdivision.} \label{BaryReflect2}
\end{center}
\end{figure}

Now, we perform some calculations. Consider two adjacent triangular tiles an any stage of barycentric subdivision. Together, they form a square. Notice in Figures \ref{BaryReflect1}, \ref{BaryReflect2} that reflection about the common edge preserves the tiling (i.e. cell structure) and swaps cuts and flows; in the language of Cannon, Floyd, and Parry, there is a weak isomorphism taking cuts to flows and vice versa \cite{subdivision}. So the optimal height and length of the square must be equal, and either combinatorial modulus of the square must be exactly one.

Now, consider the star of any vertex. It consists of some number $n$ of triangles arranged about the vertex. When we subdivide, we can remove the new star about the vertex (consisting of 2$n$ tiles) to get an annulus. This annulus consists of 2$n$ square of the type considered above, arranged with alternating orientations.

\begin{figure}
\begin{center}
\scalebox{.8}{\includegraphics{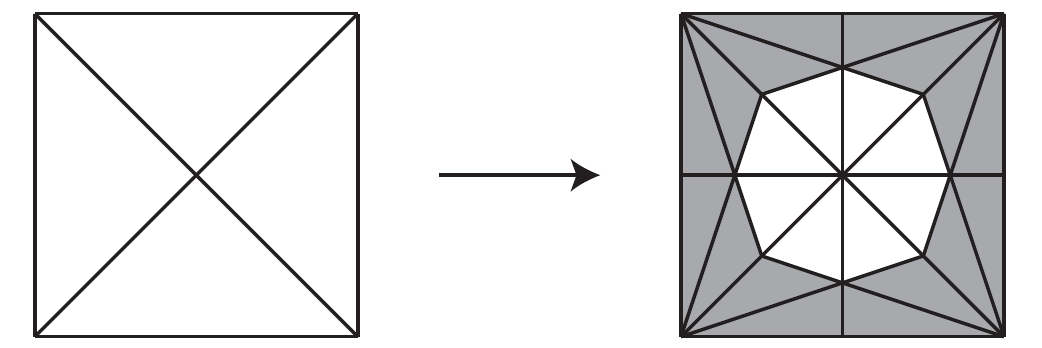}} \caption{The star of every vertex contains an annulus made of squares.} \label{BaryAnnuli}
\end{center}
\end{figure}

\begin{figure}
\begin{center}
\scalebox{.8}{\includegraphics{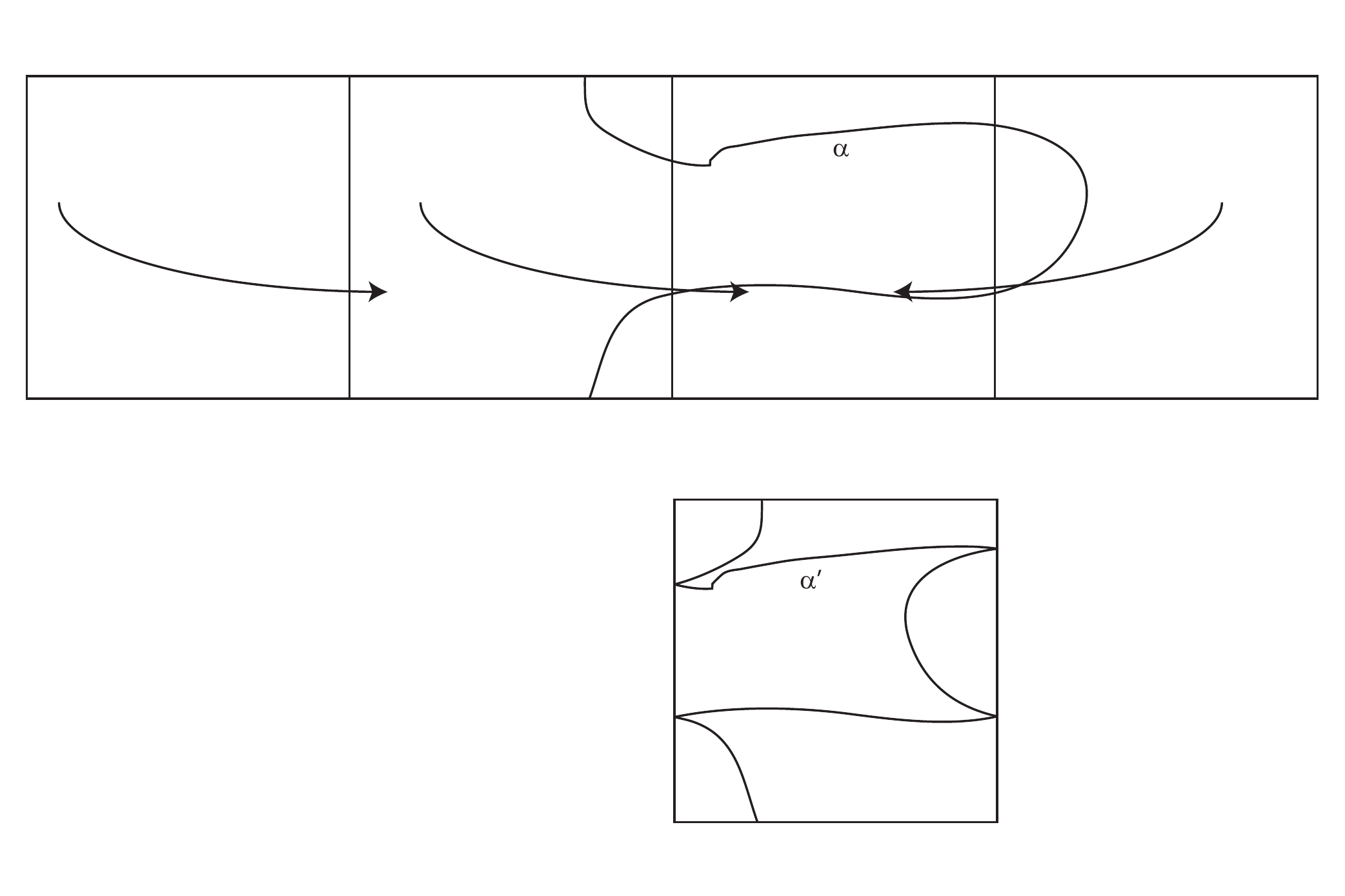}} \caption[Reflecting a cut of an annulus.]{If we consider the annulus as four squares with edges identified, we can reflect any cut of the annulus so that it lies in one given square.} \label{ReflectCut}
\end{center}
\end{figure}

Let $\alpha$ be a cut in the annulus (i.e. an arc connecting the ends of the annulus). Then pick any square in the annulus, and notice that reflecting $\alpha$ every time it touches the boundary of the square gives us an arc that remains entirely in the square, while having the exact same length (see Figure \ref{ReflectCut}). Thus, the height of the annulus is equal to the height of the square, and the are is equal to $n$ times the area of the square (after we multiply by an appropriate constant). Thus, both moduli are $\frac{1}{n}$. This class of annuli, therefore, satisfy Axiom 1.

However, we cannot layer these annuli to get larger annuli which we know to be of unbounded modulus. To see this, notice that the next smaller annulus has $2n$ squares and has modulus $\frac{1}{2n}$, and, in general, the $(k+1)$th annulus has modulus $\frac{1}{2^k n}$, and layering all of these together gives us a larger annulus which has $\mathop{\Sigma} \limits_{k=0}^\infty \frac{1}{2^k n}=\frac{2}{n}$ as a lower bound for its modulus. The exponential growth about the vertex makes the modulus shrink too quickly. This contrasts strongly with the finite valence case, in which we have nested annuli with identical moduli that sum to infinity. In the next section, we will consider subdivision rules that have linear growth at each vertex, and which thus have annuli of unbounded size.

While we don't have an upper bound for the size of the annuli, Circlepack suggests that such a limit exists. If a basis of annuli satisfy Axiom 1 and Axiom 2, the circle packings will converge to a limit function on the sphere \cite{french}. However, the packings for barycentric subdivision seem to converge to a relation in which the set of vertices is mapped to a dense set of disks. See the figure on page \pageref{CircleBary}. Note that all points that are \emph{not} vertices are contained in the double star of a new vertex at every stage of subdivision. The star of each of these vertices contains an annulus surrounding the point of modulus $\frac{1}{12}$ or $\frac{1}{8}$ for all sufficiently refined stages of subdivision. Thus, we can layer to get unbounded moduli in this case. This means that, under a sequence of circle packings, the complement of the vertices is mapped homeomorphically (by Arzela-Ascoli and the fact that locally conformal circle packings are locally equicontinuous \cite{french}) to a subset of the sphere, i.e. the complement of the disks, which is similar to a Julia set.

This is related to the work of Cannon, Floyd and Parry on subdivision rules and rational maps \cite{rational}, in which they showed that many subdivision rules can be realized as a rational map on the sphere, where the subdivision rule is obtained by pulling back a cell structure on the sphere. Conformal subdivision rules have a Julia set that is the entire sphere, but the map associated to barycentric subdivision has a Sierpsinski carpet as its Julia set.

\section{Borromean Rings}\label{BorroChap}
In contrast to barycentric subdivision, a subdivision rule associated with the Borromean rings is conformal, although the valence at each vertex remains unbounded.

\begin{figure}
\begin{center}
\scalebox{.9}{\includegraphics{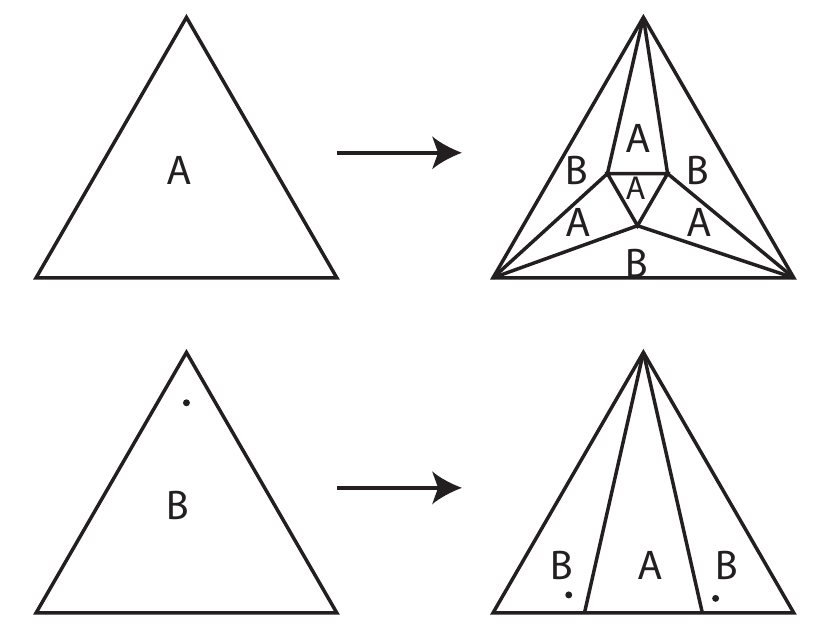}} \caption{One subdivision rule for the Borromean rings. All polyhedra are octahedra.}
\label{BorroSubs}
\end{center}
\end{figure}

The tile types are shown in Figure \ref{BorroSubs}. The dot indicates the orientation of the B tiles.

\begin{figure}
\begin{center}
\scalebox{.8}{\includegraphics{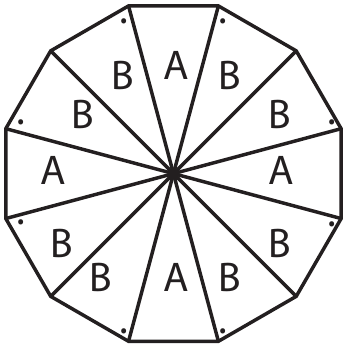}} \caption{The star of a vertex in the Borromean rings' subdivision rule.}
\label{BorroVertex}
\end{center}
\end{figure}

Let's examine the star of a vertex. Each new vertex in the subdivision rule, after all tiles around it have been subdivided once, has a star of the form shown in Figure \ref{BorroVertex}.

\begin{figure}
\begin{center}
\scalebox{.75}{\includegraphics{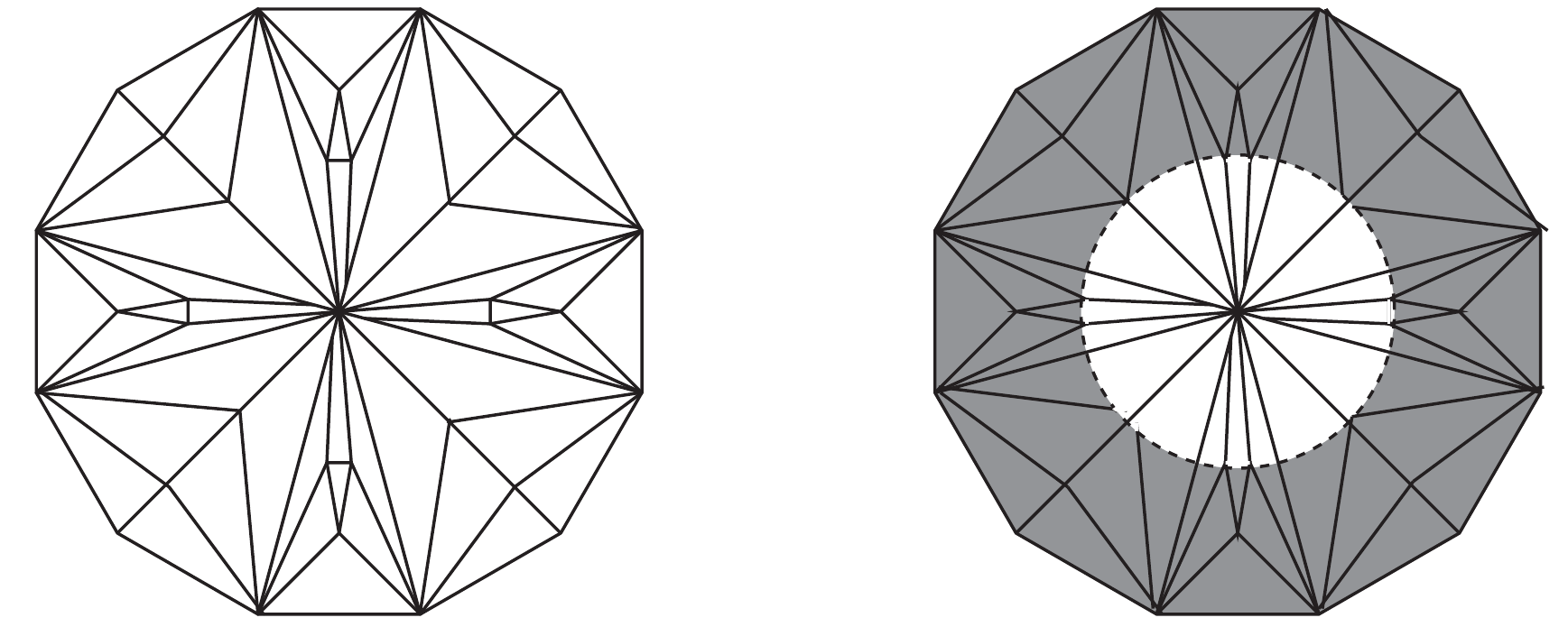}} \caption{We find an annulus in the subdivided star whose modulus is easy to calculate.}
\label{BorroVertexRings}
\end{center}
\end{figure}

After another subdivision, we have Figure \ref{BorroVertexRings}. Note the marked annulus.

\begin{figure}
\begin{center}
\scalebox{.7}{\includegraphics{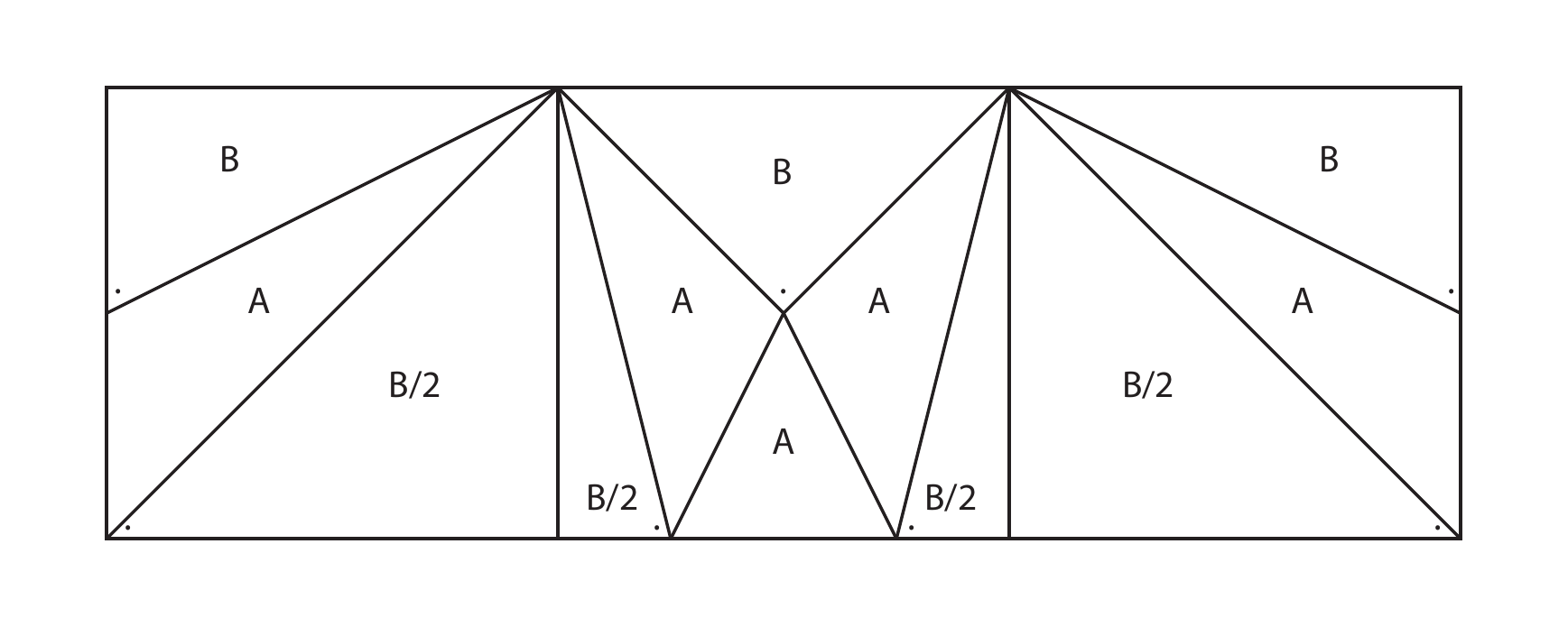}} \caption{One of the four quadrilaterals that make up the annulus in Figure \ref{BorroVertexRings}. Note that the two pairs of $B\backslash 2$ tiles don't form two whole $B$ tiles; their other halves are below the picture.}
\label{BorroQuad}
\end{center}
\end{figure}

Notice that this annulus can be reflected twice in between neighboring B tiles to get a quadrilateral of the form shown in Figure \ref{BorroQuad}. Several B tiles have had their bottom halves cut off.

\begin{figure}
\begin{center}
\scalebox{.7}{\includegraphics{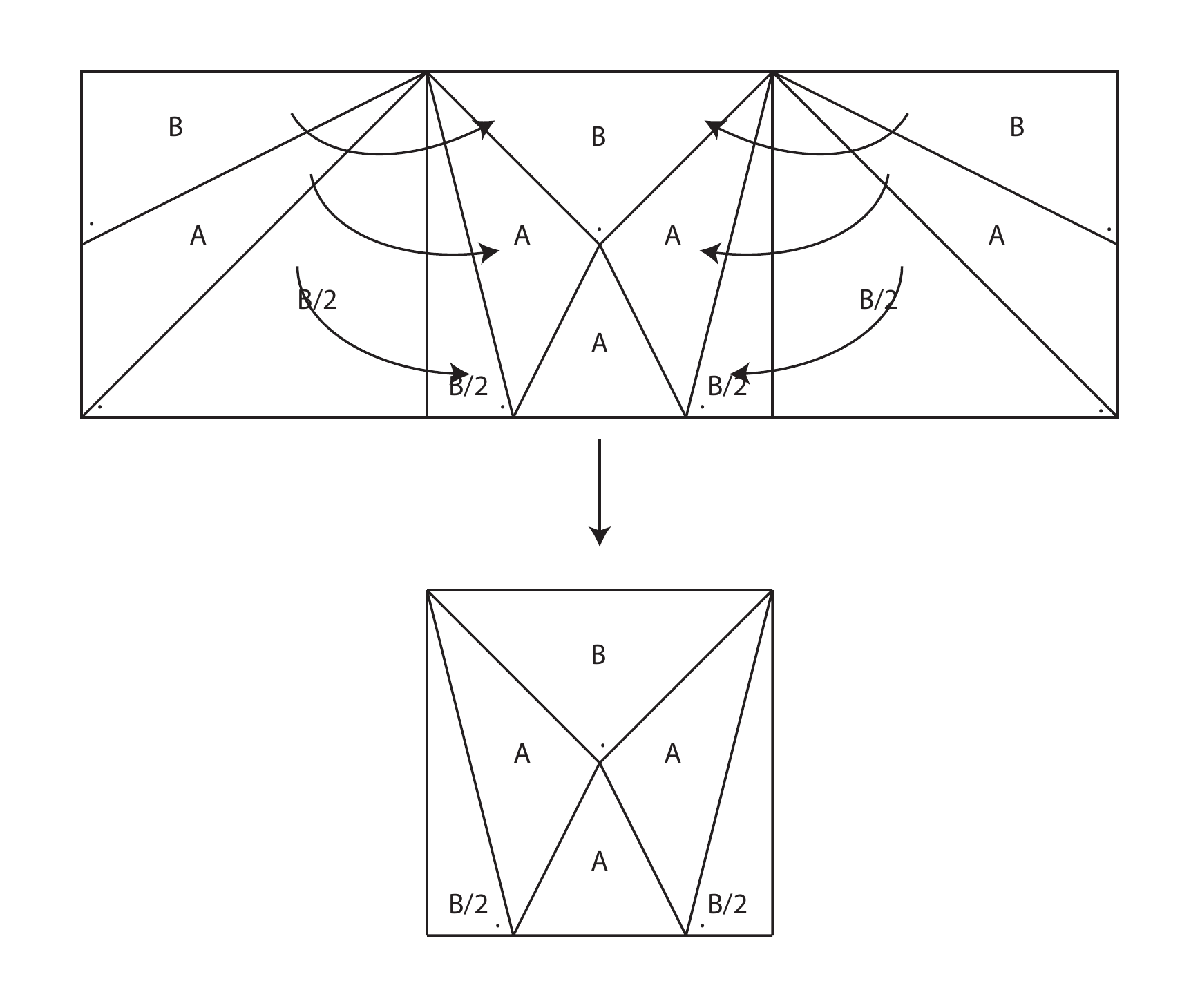}} \caption[Folding a quadrilateral.]{There is a weak cellular folding (for definition, see \cite{subdivision}) that takes the two outer tiles into inner tiles, with tops and bottoms going to tops and bottoms.}
\label{BorroReflect}
\end{center}
\end{figure}

But we can fold this up into itself as shown in Figure \ref{BorroReflect}

\begin{figure}
\begin{center}
\scalebox{.7}{\includegraphics{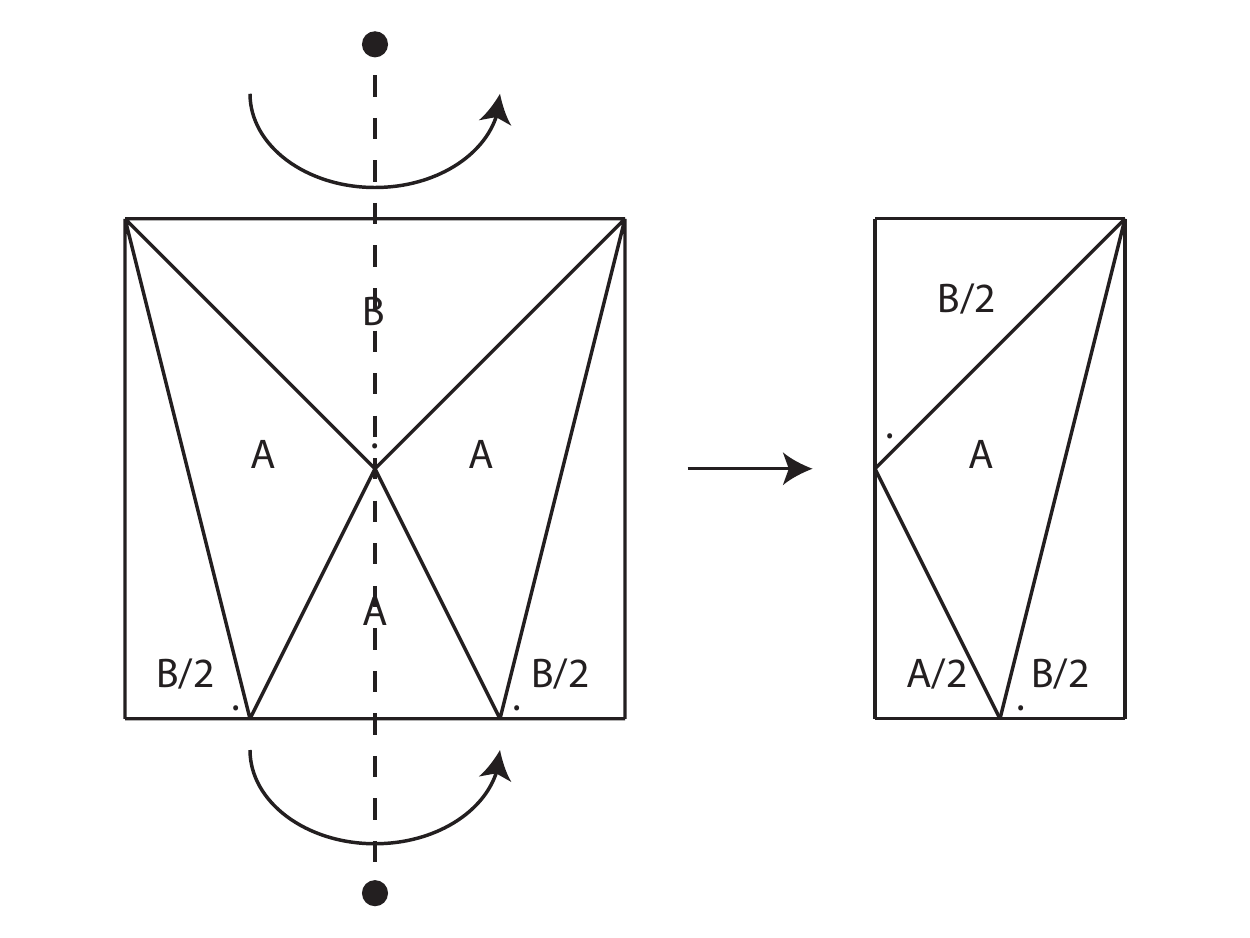}} \caption{The tile is folded along a line of symmetry.}
\label{BorroReflect2}
\end{center}
\end{figure}

This tile, in turn, can be folded in on itself, as shown in Figure \ref{BorroReflect2}.

\begin{figure}
\begin{center}
\scalebox{.7}{\includegraphics{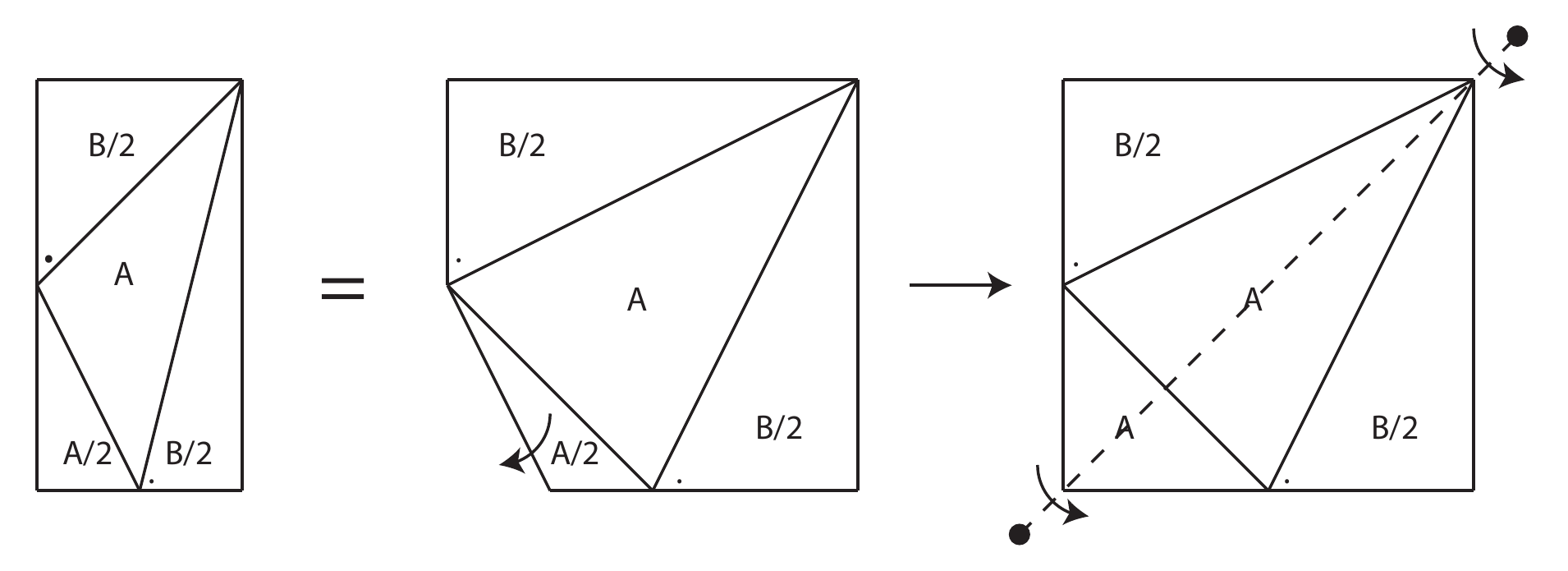}} \caption[Using reflections to calculate modulus.]{We first redraw our shape to show how it is to be unfolded, and then unfold it. The shape on the right has modulus 1, by symmetry.}
\label{BorroReflect3}
\end{center}
\end{figure}

We almost have something whose modulus we can calculate. In one final step, we unfold our shape to get something with modulus 1, as shown in Figure \ref{BorroReflect3}.

This final shape (call it Q) has modulus 1 in all stages of subdivision, as it has a symmetry interchanging the top/bottom with the left/right. We can now estimate the modulus of the original annulus by "pulling back" the weighting given to this tile. We do this by undoing each step.

Putting back the fold in Figure \ref{BorroReflect3}, we get something with the same height as Q, but with less area than that of Q. Thus, the modulus is at least 1. Undoing the steps in Figure \ref{BorroReflect2} at most doubles the area, giving a quadrilateral with modulus at least $\frac{1}{2}$. Undoing Figure \ref{BorroReflect} gives us modulus at least $\frac{1}{6}$, since area at most triples. Unfolding this into the original annulus gives modulus at least $\frac{1}{24}$.

\begin{figure}
\begin{center}
\scalebox{.45}{\includegraphics{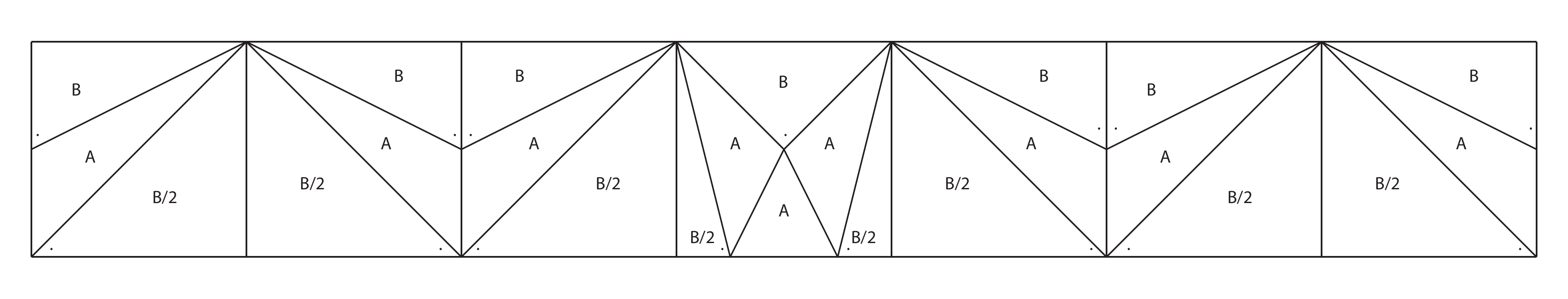}} \caption[A larger quadrilateral.]{A figure similar to Figure \ref{BorroQuad}, but for the star of a vertex at a later stage of subdivision.}
\label{BorroBigQuad}
\end{center}
\end{figure}

We can find similar annuli in the star of each vertex at any stage of subdivision, whose modulus can be calculated in the same way; the annulus at each stage can be split into four quadrilaterals which have more squares in them, as in Figure \ref{BorroBigQuad}. By reflection, we can fold in the extra tiles, and counting arguments show that the modulus is at least $\frac{1}{2n}$, where $n$ is the number of squares in the star.  The annuli thus obtained are disjoint, and so we can layer them to get annuli whose moduli are bounded below by the partial sums of a harmonic series. Assuming that the actual modulus is not much greater than this, we have logarithmic growth of modulus at each vertex. This is supported by circle packings, like the figures on pages \pageref{CircleBorro} and \pageref{CircleBorroBig}, where the stars of each vertex seem to shrink slowly.

\section{Hexagonal refinement}\label{Hexagonal}
Hexagonal refinement is a classic example of a conformal subdivision rule. The subdivision rule takes triangles and subdivides them into four smaller triangles by connecting midpoints of edges (see the figure on page \pageref{CircleHexed}). While it is easy to calculate the exact modulus of annuli in this sequence of tilings, we give lower bounds for the modulus as in our other examples.

\begin{figure}
\begin{center}
\scalebox{.7}{\includegraphics{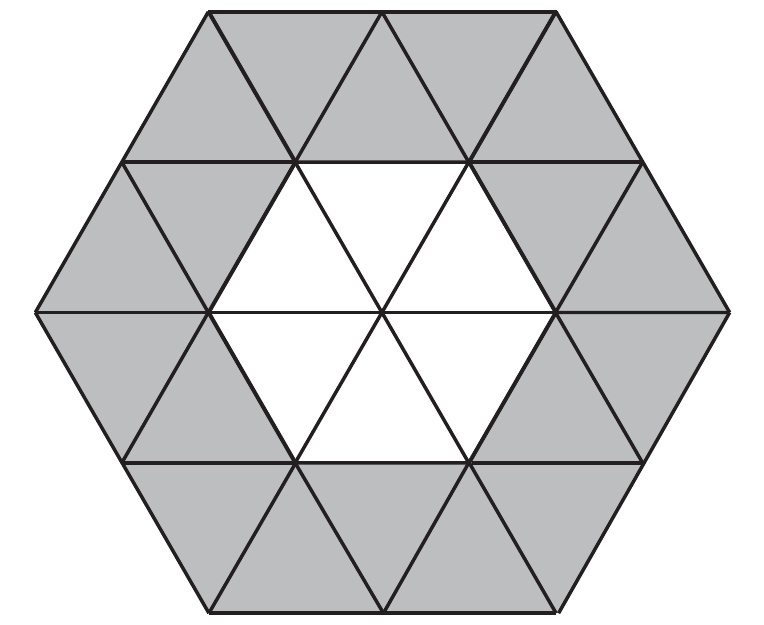}} \caption{Eventually, at every stage of subdivision, every vertex is surrounded by an annulus of this form.}
\label{HexaRing}
\end{center}
\end{figure}

\begin{figure}
\begin{center}
\scalebox{1.0}{\includegraphics{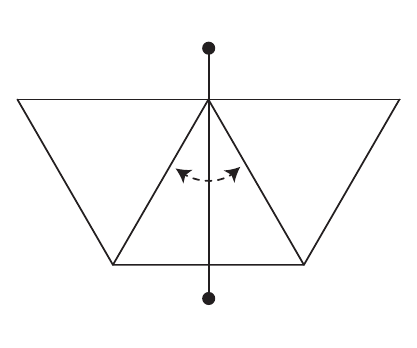}} \caption{The annulus in Figure \ref{HexaRing} is made of 6 quadrilaterals of this form. By reflection, we can assume that all cuts stay in two neighboring triangles.}
\label{HexaSpin}
\end{center}
\end{figure}

\begin{figure}
\begin{center}
\scalebox{1.4}{\includegraphics{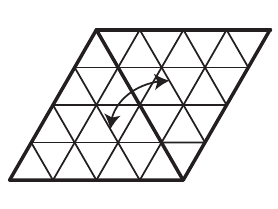}} \caption{This square has modulus one, even after several stages of subdivision.}
\label{HexaReflect}
\end{center}
\end{figure}

At every stage of subdivision, every vertex is surrounded by an annulus of 18 triangles, as shown in Figure \ref{HexaRing}. By reflection, we can consider a single quadrilateral as shown in Figure \ref{HexaSpin}. Note that, by the reflection shown in that figure, we can assume that any cuts remain in the right two triangles. These last triangles, by the reflection shown in Figure \ref{HexaReflect}, have modulus 1. Working backwards, we see that the quadrilateral in Figure \ref{HexaSpin} has modulus $\geq \frac{1}{2}$. Thus, the entire annulus in Figure \ref{HexaRing} has modulus $\geq \frac{1}{12}$.

In contrast to barycentric subdivision and the Borromean rings, the annuli surrounding a given vertex at each stage have a fixed, constant size $C\geq \frac{1}{12}$. Thus, layering gives us annuli of modulus $\geq Cn$ at stage $n$, and so axioms 1 and 2 are satisfied at each vertex. This easy layering for finite valence subdivision rules is a small part of why Cannon, Floyd, and Parry were able to replace axioms 1 and 2 with a simpler axiom 0 for finite valence subdivision rules \cite{Rich}:

\textbf{Axiom 0.} Given $x\in Y$ and a neighborhood $N$ of $x$, there is a ring $R$ in $N$ surrounding $x$ such that the moduli $m(R,T_i)$ are bounded away from 0.

Note that the estimate for every vertex in hexagonal refinement is the same at every stage; this is reflected in its circle packing on page \pageref{CircleHexed}.

\section{Generalizations}\label{Generalizations}

Notice that the key ideas in the estimates for our three subdivision rules were:
\begin{enumerate}
\item Finding pieces of modulus 1, and
\item Showing that we need only consider cuts that remain in these pieces.
\end{enumerate}

One powerful technique related to these ideas is the 1,2,3-tile criterion of Cannon, Floyd, and Parry \cite{subdivision}. Essentially, a subdivision rule satisfies the criterion if every `test quadrilateral' formed from one or two tiles in a row (or three if the middle one is a triangle) has modulus uniformly bounded away from 0. This criterion allows us to argue as we did in the three cases above to provide estimates.

There are three types of test quadrilaterals:

\emph{Type I}. The quadrilateral $Q$ consists of one tile, and each of the ends of $Q$ consists of one edge.

\emph{Type II}. The quadrilateral $Q$ consists of two tiles whose intersection consists of one edge $f$, and each of the ends of $Q$ consists of one edge which meets $f$. We call $f$ the \textbf{interior edge} of $Q$.

\emph{Type III}. The quadrilateral $Q$ consists of three tiles $t_1,t_2$ and $t_3$, where $t_2$ is a triangle with edges $f_1$, $f_2$, and $f_3$. The intersection of $t_1$ and $t_2$ is $f_1$, and the intersection of $t_2$ and $t_3$ is $f_3$. The intersection of $t_1$ and $t_3$ is a vertex $v$. The top of $Q$ consists of an edge containing $v$, and the bottom of $Q$ is $f_2$. We call $f_1$ and $f_3$ the \textbf{interior edges} of $Q$.

Cannon, Floyd and Parry showed that a uniform lower bound on the modulus of all test quadrilaterals (call it M) gives a lower bound for the modulus of a ring R in our surface, depending on the size of the ring. Specifically, in the proof of Theorem 5.1 in \cite{subdivision}, they take the star of a simple closed curve $\alpha$ in some subdivision $R^i$ of the ring, where $\alpha$ misses every vertex, and show that any flow curve from the bottom of this star to the top must join the ends of a test quadrilateral (i.e. its image in some test quadrilateral is a cut). They then give the star of $\alpha$ a weight vector that is the sum of all optimal weightings of test quadrilaterals, normalized to have area 1. The star of $\alpha$, then, has modulus $\geq \frac{M}{Ak}$, where $k$ is the number of tiles in the star of $\alpha$ in $R^i$ and $A$ is the maximum area a single tile can have under the weighted sum (which is seen to be bounded, as the number of test quadrilaterals a single tile can lie in is uniformly bounded).

When the star of a vertex $v$ in a subdivided surface $S$ is a disk, we can take its boundary and push it slightly outward to be such an $\alpha$. The number of tiles in the star of $\alpha$ will then be $\leq B_i(v) val_i(v)$, where $B_i(v)$ is the maximum valence of vertices surrounding $v$. In a subdivision rule with combinatorial mesh going to 0, each vertex $v$ has a uniform bound on $B_i(v)$, as each edge surrounding the vertex is subdivided periodically. Thus, we see that any annulus containing the star of $\alpha$ must have modulus $\geq \frac{C}{val_i(v)}$ for some constant $C$..

By nesting this sequence of annuli, we can construct an annulus of larger modulus. When valence grows exponentially at each vertex, as in barycentric subdivision, we get a geometric sequence of estimates of the size of annuli, and we cannot construct annuli whose modulus estimates are unbounded. I conjecture that every subdivision rule with exponential growth at each vertex which satisfies the 1,2,3-tile criterion will satisfy axiom 1 but not axiom 2, just like barycentric subdivision. In this case, the limit of circle packings will be a continuous relation instead of a continuous function onto the sphere, where points are blown up into disks. In \cite{myself2}, there are subdivision rules for hyperbolic 3-manifolds with hyperbolic surface boundary; I conjecture that all such subdivision rules satisfy the 1,2,3-tile criterion and behave like barycentric subdivision at each vertex corresponding to a hyperbolic surface in the boundary.

Nesting annuli in a subdivision rule of linear growth gives us a harmonic series for our estimates, which shows that all such subdivision rules that satisfy the 1,2,3-tile criterion are conformal. I conjecture that all the subdivision rules arising from alternating links in \cite{linksubs} satisfy this criterion, and are conformal.

Finally, subdivision rules of bounded valence give a linear modulus estimate, and those that satisfy the 1,2,3-tile criterion are conformal, which was why Cannon, Floyd and Parry originally proved the criterion in \cite{subdivision}.

\section{Circle packed pictures}\label{CirclePacked}
All the pictures on pages \pageref{CircleHexed}-\pageref{CircleBorroBig} were created with Ken Stephenson's Circlepack \cite{Circlepak}.

\begin{figure}
\begin{center}
\scalebox{0.8}{\includegraphics{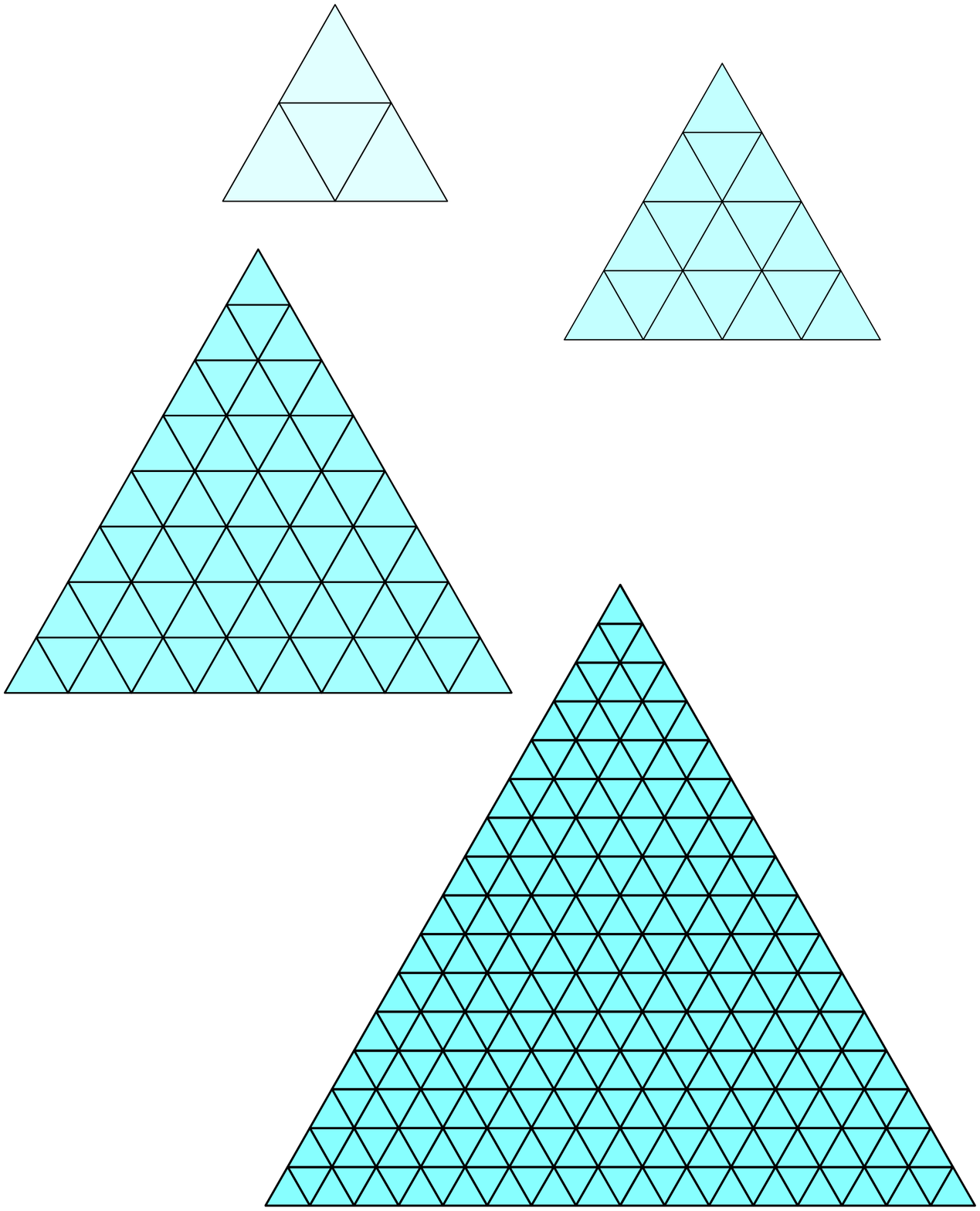}}
\label{CircleHexed}
\end{center}
\end{figure}

\begin{figure}
\begin{center}
\scalebox{0.8}{\includegraphics{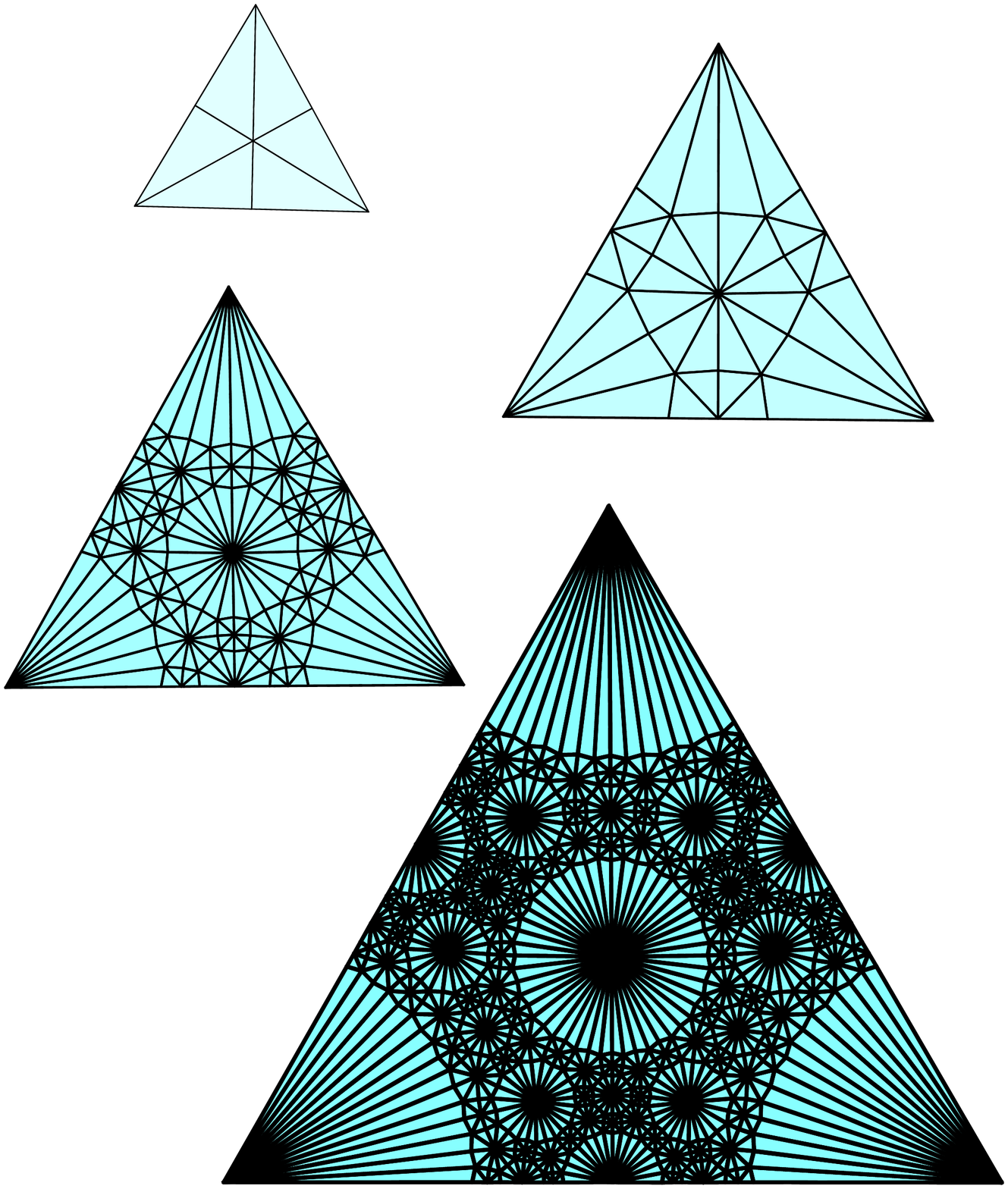}}
\label{CircleBary}
\end{center}
\end{figure}

\begin{figure}
\begin{center}
\scalebox{0.8}{\includegraphics{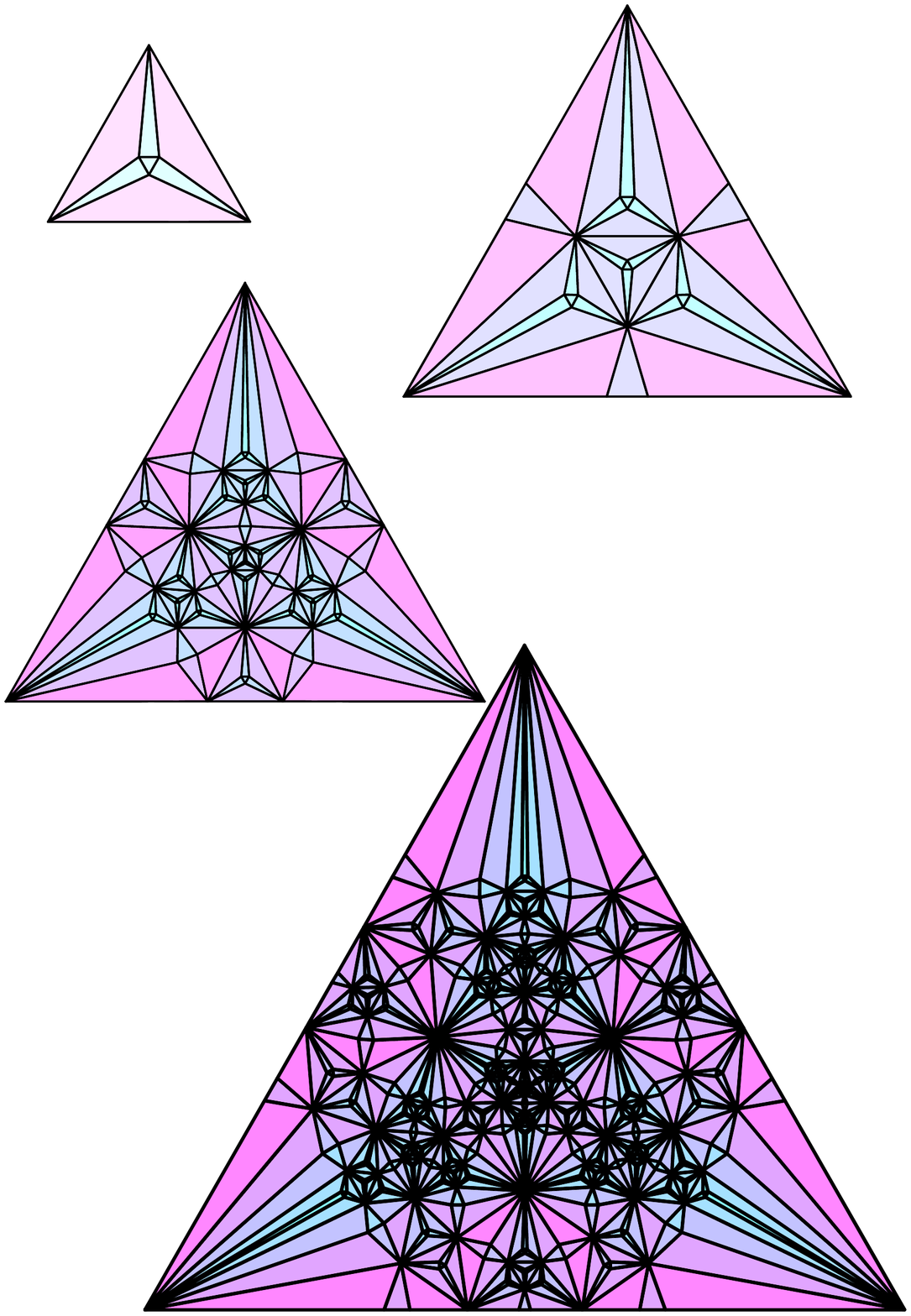}}
\label{CircleBorro}
\end{center}
\end{figure}

\begin{figure}
\begin{center}
\scalebox{0.8}{\includegraphics{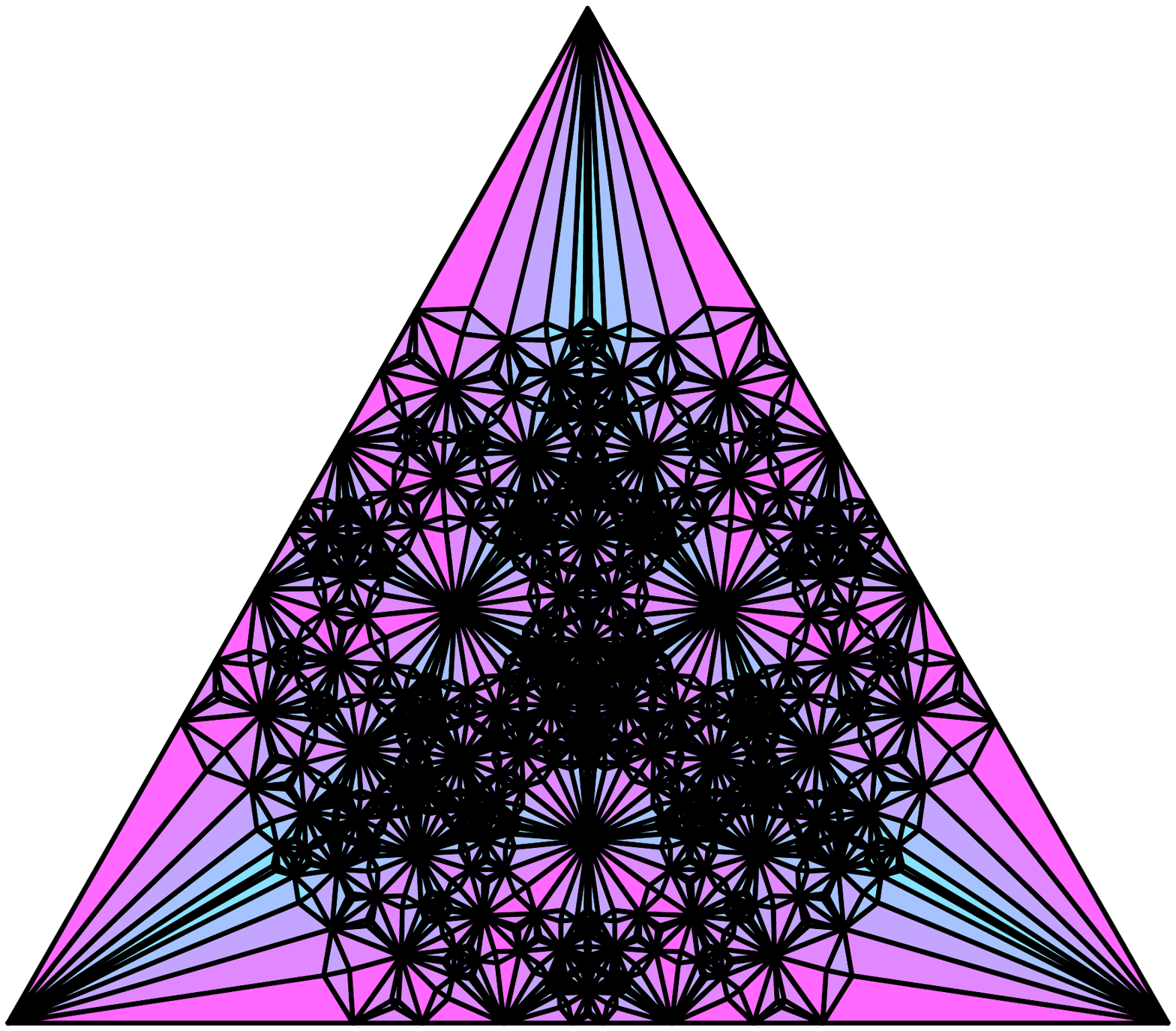}}
\label{CircleBorroBig}
\end{center}
\end{figure}

\bibliographystyle{plain}
\bibliography{ValPaper}

\end{document}